\documentclass[12pt,english,british]{article}
\usepackage[T1]{fontenc}
\usepackage[latin9]{inputenc}
\usepackage{geometry}
\usepackage{amsmath}
\usepackage{amsthm}
\usepackage{amssymb, verbatim}

\usepackage{hyperref}

\makeatletter
\theoremstyle{plain}
\newtheorem{thm}{\protect\theoremname}[section]
\newtheorem{cor}[thm]{Corollary}

  \theoremstyle{definition}
  
  \theoremstyle{definition}
  \newtheorem{example}[thm]{\protect\examplename}
  \theoremstyle{remark}
  \newtheorem{rem}[thm]{\protect\remarkname}

\newcommand{\N}{\mathbb{N}}
\newcommand{\Z}{\mathbb{Z}}
\newcommand{\R}{\mathbb{R}}
\newcommand{\FF}{\mathcal{F}}
\newcommand{\ff}{\mathbf{f}}
\newcommand{\TD}{\dim_T}
\newcommand{\del}{\partial}
\newcommand{\butnot}{\setminus}
\newcommand{\hh}{\mathbf{h}}
\newcommand{\PP}{\mathcal P}

\newcommand{\pp}{\mathbf{p}}
\newcommand{\HD}{{\dim_H}}

\newcommand{\BD}{{\dim_B}}
\newcommand{\xx}{\mathbf{x}}
\newcommand{\bx}{\mathbf{x}}
\newcommand{\yy}{\mathbf{y}}
\newcommand{\AD}{\dim_A}
\renewcommand{\gg}{\mathbf{g}}
\newcommand{\eps}{\varepsilon}
\newcommand{\sep}{\,:\,}

\newcommand{\bfm}{\mathbf{m}}

\global\long\def\subp{\supseteq}
\global\long\def\bbr{\mathbb{R}}
\global\long\def\bbz{\mathbb{Z}}

\allowdisplaybreaks
\makeatother

\usepackage{babel}
  \addto\captionsbritish{\renewcommand{\definitionname}{Definition}}
  \addto\captionsbritish{\renewcommand{\examplename}{Example}}
  \addto\captionsbritish{\renewcommand{\remarkname}{Remark}}
  \addto\captionsbritish{\renewcommand{\theoremname}{Theorem}}
  \addto\captionsenglish{\renewcommand{\definitionname}{Definition}}
  \addto\captionsenglish{\renewcommand{\examplename}{Example}}
  \addto\captionsenglish{\renewcommand{\remarkname}{Remark}}
  \addto\captionsenglish{\renewcommand{\theoremname}{Theorem}}
  \providecommand{\definitionname}{Definition}
  \providecommand{\examplename}{Example}
  \providecommand{\remarkname}{Remark}
\providecommand{\theoremname}{Theorem}

\RequirePackage[usenames,dvipsnames]{xcolor}

\newif\ifdraft\drafttrue

\newif\ifcolorcomments

\colorcommentstrue

\newcommand{\allowcomments}[4]{
\newcommand{#1}[1]{\ifdraft{\ifcolorcomments{ \textcolor{#4}{##1 --#3}}\else{\textsl{ ##1 \ --#3}}\fi}\else{}\fi}
}

\allowcomments{\comabel}{AF}{Abel}{blue}
\allowcomments{\comjon}{JF}{Jonathan}{red}
\allowcomments{\comerez}{EN}{Erez}{blue}
\allowcomments{\comdavid}{DS}{David}{Green}

\title{Schmidt's game on Hausdorff metric and function spaces: generic dimension of sets and images}

\author{\'Abel Farkas$^a$, Jonathan M. Fraser$^b$, Erez Nesharim$^c$, \& David Simmons$^d$}

\begin{document}

\maketitle

\vspace{-8mm}
\begin{center}
	$^a$Alfr\'ed R\'enyi Institute of Mathematics,  Hungarian Academy of Sciences, Budapest, Hungary.\\
	 \vspace{3mm}
	$^b$Mathematical Institute, University of St Andrews, UK.\\
	 \vspace{3mm}
    $^c$Faculty of Mathematics, Technion, Israel.\\
     \vspace{3mm}
    $^d$Department of Mathematics, York, UK.\\
     \vspace{3mm}
	\end{center}
\begin{abstract}
We consider Schmidt's game on the space of compact subsets of a given metric space equipped with the Hausdorff metric, and the space of continuous functions equipped with the supremum norm. We are interested in determining the generic behaviour of objects in a metric space, mostly in the context of fractal dimensions, and the notion of `generic' we adopt is that of being winning for Schmidt's game. We find properties whose corresponding sets are winning for Schmidt's game that are starkly different from previously established, and well-known, properties which are generic in other contexts, such as being residual or of full measure.

\emph{Mathematics Subject Classification} 2010:  primary: 28A80, 91A44; secondary:   28A78, 91A05.

\emph{Key words and phrases}: Schmidt's game, dimension.
\end{abstract}

\section{Schmidt's game and winning sets}

We consider Schmidt's game introduced by Schmidt in \cite{wmS66}. 
The game is played in a complete metric space $\left(X,d\right)$ and it has some similarities to the Banach--Mazur game. The game is played by two players, Alice and Bob, and the rules are described below. Given $0<\alpha,\beta<1$, Alice
and Bob play the \emph{$\left(\alpha,\beta\right)$-game} as follows:
\begin{enumerate}
\item Bob begins by choosing $r_0>0$ and $x_0\in X$. Write $B_0=B\left(x_0,r_0\right)$ for the closed ball of radius $r_0$ centered at $x_0$
and $r_n=\left(\alpha\beta\right)^nr_0$ for all $n\in\mathbb{\mathbb{N}}$.
\item On Alice's $n$th turn ($n \geq 1$), she chooses $y_n$ such that $d\left(x_{n-1},y_n\right)+\alpha r_{n-1}\leq r_{n-1}$.
\item On Bob's $n$th turn ($n \geq 1$), he chooses $x_{n}$ such that
$d\left(y_n,x_{n}\right)+\alpha\beta r_{n-1}\leq\alpha r_{n-1}$.
\item The inequalities above ensure that the closed balls $\left\{ B_n=B\left(x_n,r_n\right) \right\}_{n=0}^\infty $ and $\left\{ A_n=B\left(y_n,\alpha r_{n-1}\right) \right\}_{n=1}^\infty $
form a decreasing sequence that satisfies
\begin{equation}\label{eq:seqOfBalls}
B_0\subp A_1\subp B_1 \subp A_2\subp\cdots\subp B_n \subp A_{n+1} \subp B_{n+1}\subp \cdots,
\end{equation}
and hence intersect at a unique point which is called the \emph{outcome} of the game.
\end{enumerate}
Given a set $S \subseteq X$, if Alice has a strategy guaranteeing that the outcome
lies in $S$, then $S$ is called \emph{$\left(\alpha,\beta\right)$-winning}.
If for some fixed $\alpha$, the set $S$ is $\left(\alpha,\beta\right)$-winning for all $0<\beta<1$, then $S$ is called \emph{$\alpha$-winning}.
If $S$ is $\alpha$-winning for some $0<\alpha<1$, then $S$ is called
 \emph{winning}. Sets which are winning should be thought of as being `big' (cf. Theorem \ref{thm:propertiesOfWinning} below). A property is \emph{generic with respect to Schmidt's game} if the set of points with that property is winning.  The Banach--Mazur game is played similarly, but without any restriction on the size of the nested sequence of balls \eqref{eq:seqOfBalls}, that is, the Banach--Mazur game is a topological game, whereas Schmidt's game is metric.
A fundamental fact which is important to keep in mind throughout this paper is that the winning sets for the Banach--Mazur game are precisely the residual sets, that is, sets whose complement is a countable union of nowhere dense sets.  We refer the reader to \cite{ox} for a precise definition of the Banach--Mazur game and more on its properties.

A common application of the Banach--Mazur game is that if a countable collection of properties are winning for the Banach--Mazur game then it is also winning that all of them are satisfied at once. Schmidt's game can be applied similarly as the following theorem shows. We list some well-known properties of Schmidt's game. For reference see \cite{wmS66}.

\begin{thm}\label{thm:propertiesOfWinning}
Winning sets have the following
properties:
\begin{enumerate}
\item If $S\subseteq X$ is winning, then $S$ is dense. If $S$ is $\alpha$-winning
for $\alpha>\frac12$, then $S=X$.
\item If $S\subseteq X$ is $\alpha$-winning and $0<\alpha'\leq\alpha$,
then $S$ is $\alpha'$-winning.
\item If $S\subseteq\mathbb{R}^{d}$ is winning, then the Hausdorff dimension of $S$ is $d$.
\item If $S_k\subseteq X$ is $\alpha$-winning for every $k\in\mathbb{N}$,
then $S=\bigcap_{k=1}^{\infty}S_k$ is $\alpha$-winning.
\item If $S_1,\ldots,S_N\subseteq X$ are winning sets ($N$ finite) then $S=\bigcap_{k=1}^{N}S_k$ is winning.
\item If $S \subseteq X$ is winning, then $X \setminus S$ is not winning.
\end{enumerate}
\end{thm}
The property of being winning is subtle and how it relates to other notions of being generic is of particular interest and a key theme of this paper. For example, the set of badly approximable numbers
\[
\mathrm{BA}=\left\{ x\in\mathbb{R} \sep \inf_{q\in\mathbb{N},\,p\in\Z} \, q \,   | qx - p |>0\right\}
\]
is well-known to be of first category (co-residual) and have zero Lebesgue measure,
but nevertheless be winning. Also, the set of numbers which are not normal to some fixed base is well-known
to have zero Lebesgue measure, and nevertheless to be residual and winning, see \cite{borel, wmS66, volkmann}.

\section{Results}

\subsection{Dimensions of sets in a metric space}

Let $(X,d)$ be a complete separable metric space. Let $\mathcal{K}(X)$ be the set of nonempty compact subsets of $X$ equipped with the
Hausdorff metric, denoted by $d_H$. It is well-known that $(\mathcal{K}(X),d_H)$ is complete. Schmidt's game will be played on $\mathcal{K}(X)$ and the winning properties of subsets of $\mathcal{K}(X)$ defined according to dimension will be considered. Five standard notions of dimension will be used: the lower, Hausdorff, lower box, upper box, and Assouad dimensions, denoted by $\dim_{L}$, $\dim_{H}$, $\underline{\dim_B}$, $\overline{\dim_B}$, and $\dim_{A}$,  respectively.  We note that the lower dimension is sometimes referred to as the lower Assouad dimension in the literature.


We recall the definitions here for convenience.  Fix a non-empty subset $K \subseteq X$ and write $\mathrm{diam}(A) \in [0,\infty]$ for the diameter of a set $A$.  Given $s \geq 0$, the  \emph{$s$-dimensional Hausdorff (outer) measure} of $K$ is defined by
\[
\mathcal{H}^s (K) = \lim_{ r \to 0}\inf \left\{ \sum_{i } {\mathrm{diam}(U_i)}^s : \{ U_i \}_{i} \text{ is a countable $ r$-cover of $K$} \right\},
\]
where the $r$-cover of a set $K$ is a family of sets with diameter at most $r$ such that the family is a cover of $K$. The \emph{Hausdorff dimension} of $K$ is then
\[
\dim_{H}  K = \inf \left\{ s \geq 0: \mathcal{H}^s (K) =0 \right\} = \sup \left\{ s \geq 0: \mathcal{H}^s (K) = \infty \right\}.
\]
The \emph{lower and upper box dimensions} of $K$ are defined by
\[
\underline{\dim_B}  K = \liminf_{ r \to 0} \, \frac{\log N_ r (K)}{-\log  r} \qquad
\text{and}
\qquad
\overline{\dim_B} K = \limsup_{ r \to 0} \,  \frac{\log N_ r (K)}{-\log  r}\,,
\]
respectively, where $N_ r (K)$   denotes the \emph{minimal number of open sets of diameter at most $r$ required to cover $K$}.  If $\underline{\dim}_\text{B} K = \overline{\dim}_\text{B} K$, then we call the common value the \emph{box dimension} of $K$ and denote it by $\dim_\text{B} K$. The \emph{Assouad dimension} of $K$ is defined by 
\begin{eqnarray*}
\dim_\text{A} K =  \inf \Bigg\{ \  \alpha &:& \text{     there exists a constant $C >0$ such that,} \\
&\,& \hspace{20mm}  \text{for all $0<r<R $ and $x \in K$, } \\ 
&\,&\hspace{18mm}  \text{$ N_r\left( B(x,R) \cap K \right) \ \leq \ C \left(\frac{R}{r}\right)^\alpha$ } \Bigg\}\,.
\end{eqnarray*}

The \emph{lower dimension} of $K$   is defined by
\begin{eqnarray*} 
\dim_{L} K = \sup \Bigg\{ \  \alpha &:& \text{     there exists a constant $C >0$ such that,} \\
&\,& \hspace{0mm}  \text{for all $0<r<R \leq \max\{ |K|,1\} $ and $x \in K$, } \\ 
&\,&\hspace{23mm}  \text{$ N_r\left( B(x,R) \cap K \right) \ \geq \ C \left(\frac{R}{r}\right)^\alpha$ } \Bigg\}\,.
\end{eqnarray*}


The reader is referred to \cite{falconer, robinson, kaenmaki} for more discussion  of these dimensions and their basic properties. It is of particular importance in the subsequent analysis that for any compact set $K$ in a metric space the following inequalities hold:
\[
\dim_{L}K \leq \dim_{H}K \leq \underline{\dim_B}K \leq \overline{\dim_B}K \leq \dim_{A} K\,.
\]
Although these inequalities are straightforward to establish, we refer the reader to \cite{larman} for the first inequality, \cite{falconer} for the middle two, and \cite{robinson} for the final inquality.

Recall that a metric space is \emph{doubling} if and only if $\dim_{A}X < \infty$, see \cite[Lemma 9.4]{robinson}, and \emph{uniformly perfect} if and only if $\dim_{L}X > 0$, see \cite[Lemma 2.1]{kaenmaki}. We say that $X$ is \emph{Assouad sharp} if it is doubling and there exists $C>0$ such that all $x \in X$ and $0<r<R$ satisfy
\[
N_r\left( B(x,R) \right) \leq C \left( \frac{R}{r} \right)^{\dim_{A}X}\,.
\]
Examples of Assouad sharp spaces include Ahlfors regular metric spaces, in particular $\mathbb{R}^d$. Note that for Ahlfors regular spaces all of the above dimensions coincide, see \cite{bylund}.
\begin{rem}
If $\left\{B\left(x_{i},r\right)\right\}_{i=1}^{N}$ is a maximal collection of disjoint balls of
radius $r$ with centers in $K$ then $K\subseteq\bigcup_{i=1}^{N}B\left(x_{i},2r\right)$.
On the other hand any cover of $K$ with balls of radius $r$ must
contain at least $N$ balls. Hence if $\dim_{A}(K)<\infty$ then the
minimum number of balls of radius $r$ needed to cover $K$ and the
maximum number of disjoint balls of radius $r$ centered in $K$ are comparable
up to constant factor. This constant factor could be included in the constant $C$ above.
Throughout the paper we might switch back and forth between disjoint
balls and covering balls, which is a common strategy in dimension theory.
\end{rem}
\begin{thm} \label{main} \hspace{1mm}
\begin{enumerate}
\item[\textup{(i)}] If $X$ is uniformly perfect, then the set $\left\{ K\in\mathcal{K}(X)\sep\dim_{L}K>0\right\} $ is winning.

\item[\textup{(ii)}] For all $\varepsilon>0$, the set $\left\{ K\in\mathcal{K}(X)\sep\dim_{L}K<\dim_{L}X-\varepsilon\right\} $ is not winning.

\item[\textup{(iii)}] If $X$ is Assouad sharp, then the set $\left\{ K\in\mathcal{K}(X)\sep\dim_{A}K<\dim_{A}X\right\}$ is winning.

\item[\textup{(iv)}] If $X$ is doubling, then  for all $\varepsilon>0$ the set $\left\{ K\in\mathcal{K}(X)\sep\dim_{A}K >\varepsilon\right\} $ is not winning.

\end{enumerate}
\end{thm}

\begin{rem}\label{abint}
Assume that $0<a<b<c=\dim_{L}X=\dim_{A}X<\infty$. Without proof we note that Bob, while choosing $\beta$ small enough depending on $\alpha$, can use a modification of his strategies in the proof of Theorem \ref{main} (ii) and (iv) to ensure that $a<\dim_{L}K\leq \dim_{A}K<b$ for the outcome of the game, $K$. This implies that if for sets $A,B\subseteq (0,c)$, the set
\begin{equation}\label{eq:question1}
\left\{ K\in\mathcal{K}(X)\sep\dim_LK\in A,\, \dim_A K\in B,\, \dim_{L}K<\dim_{A}K\right\}
\end{equation}
is winning, then both $A$ and $B$ are dense in the interval $(a,b)$. We don't know if the converse holds, nor if there exists any other characterisation of the winning property of the set in \eqref{eq:question1} in terms of $A$ and $B$ both being `big' subsets of $(0,c)$, where `big' is chosen appropriately, for example `big' might mean full Lebesgue measure, winning, co-meager, or uncountable. 
This is not clear even for $X=\mathbb{R}^d$.
It is even less clear in the case when $X$ is uniformly perfect and doubling but satisfies $\dim_{L}X<\dim_{A}X$. We also do not know whether the property $\dim_{L}K< \dim_{A}K$ is winning even in the $X=\mathbb{R}^d$ case. However, we suspect that once $\alpha$ and $\beta$ are fixed then Alice has enough impact to force the lower and Assuad dimensions to be distinct for the outcome $K$.
\end{rem}

Part (iii) of  Theorem \ref{main} has  the weakness that we assumed $X$ to be Assouad sharp which is often not the case. It is possible to drop the Assouad sharpness by assuming instead that $X$ is uniformly perfect, though the proof becomes more complicated.

\begin{thm}\label{extension}
Let $X$ be a doubling uniformly perfect space. Then:
\begin{enumerate}
\item[\textup{(i)}] $\left\{ K\in\mathcal{K}\left(X\right) \sep \AD K<\AD X\right\}$ is winning.

\item[\textup{(ii)}] If furthermore $\overline{\dim_B} X<\infty$, then $\left\{ K\in\mathcal{K}\left(X\right) \sep \underline{\dim_B} K<\overline{\dim_B} X\right\}$ is winning.
\end{enumerate}
\end{thm}

Part (i) of Theorem \ref{main} is sharp in the sense that the statement is not true if the `uniformly perfect' assumption is dropped. For example, if $X$ has an isolated point, then the collection of sets with positive dimension cannot be winning because Bob can choose the initial set to be the singleton consisting of the isolated point and the initial radius to be so small that the outcome of the game is already forced to be the isolated point. However, one can construct a perfect set which also behaves similarly.

\begin{example}
Let $X$ be the disjoint union of the unit ball in $\mathbb{R}^d$ and a perfect compact set $A$ of Hausdorff dimension $0$ (or even of Assouad dimension  $0$). Then Bob can choose the initial set to be contained in $A$ and choose $r_0$ to be small enough to make sure that every step stays inside $A$. Hence the outcome of the game is in $A$ and so of $0$ dimension. On the other hand, having $0$ dimension is not winning because Bob can also make sure that the game is played only inside the unit ball after the initial step (where Theorem \ref{main} applies). The example shows that assuming that, for example, $\dim_H X>0$ is not enough to conclude that $\dim_H K>0$ is a winning property.
\end{example}

The converse of part (i) of Theorem \ref{main} is not true, as can be seen from the following example of a set which is not uniformly perfect, but yet the conclusion of part (i) holds.

\begin{example} \label{eg1}
Let
\[
X \ = \  \{ 0 \} \cup  \bigcup_{n=1}^\infty \  \left[1/n-2^{-n}, \ 1/n+2^{-n}\right]
\]
equipped with the Euclidean distance.  It follows from \cite[Example 2.5]{fraser} that $\dim_\text{L} X = 0$ and so $X$ is not uniformly perfect.  However, no matter which initial set $K_0 \subseteq X$  and radius $r_0$ Bob picks, as long as $\alpha < 1/2$, eventually the radius $r_n=\left(\alpha\beta\right)^nr_0$ will be small enough such that Alice can choose a set $K_n'$  such that elements of $B_H\left(K_n', \alpha r_n\right)$ are uniformly bounded away from 0.  Thus the problem reduces to the case when $X=[0,1]$ (more accurately,  a finite collection of intervals) and it follows from Theorem \ref{main} (i) that $\left\{ K\in\mathcal{K}(X)\sep\dim_{L}K>0\right\} $ is winning.
\end{example}

Part (i) of Theorem \ref{main} is also sharp in another sense, since part (iv) shows that (in the doubling case) we cannot replace $\dim_{L}K>0$ with $\dim_{L}K>\varepsilon$ for any positive $\varepsilon>0$.

Part (iii) of Theorem \ref{main} is sharp, since part (ii) shows that (provided $\dim_{L}X = \dim_{A}X$, e.g. if $X$ is Ahlfors regular) the assumption $\dim_{A}K<\dim_{A}X$ cannot be replaced with $\dim_{A}K<\dim_{A}X-\varepsilon$ for any positive $\varepsilon>0$.

To emphasise an important setting where our results are complete, we state the specialisation of Theorem \ref{main} to Euclidean space, noting that the same results also hold in any Ahlfors regular space, where $d$ is replaced by the Hausdorff dimension of the space.

\begin{cor} \label{corcor}
The set
\[
\left\{ K\in\mathcal{K}\left(\mathbb{R}^d\right)\sep0<\dim_{L}K \leq \dim_{A} K < d\right\}
\]
is winning, but for all $\varepsilon>0$ the sets  $\left\{ K\in\mathcal{K}\left(\mathbb{R}^d\right)\sep \dim_{A} K >\varepsilon\right\} $ and
$\left\{ K\in\mathcal{K}\left(\mathbb{R}^d\right)\sep\dim_{L} K<d-\varepsilon\right\}$ are not winning.
\end{cor}

This should be compared with the well-known result of Feng and Wu \cite{fengwu}, which considers the Banach--Mazur game instead of Schmidt's game. One can see that the results are rather different.

\begin{thm}[Feng--Wu, 1997]
The set
\[
\left\{ K\in\mathcal{K}\left(\mathbb{R}^d\right)\sep \dim_{L} K =\underline{\dim_{B}} K = 0 \text{ and } \overline{\dim_{B}} K =\dim_{A} K= d \right\}
\]
is winning for the Banach--Mazur game.
\end{thm}

\begin{rem}
McMullen \cite{M} introduced the concept of \emph{absolute winning} in $\mathbb{R}^n$, which can be extended to any complete metric space (see, e.g., \cite{BHNS}). In the
absolute
game
in $(X,d)$,
Bob chooses an initial closed ball $B_{0}$
of radius $r_{0}>0$ and
$\beta\in(0,1)$.
Alice then chooses a closed ball $A_{1}$ of radius at most $\beta r_{0}$, Bob chooses, if possible,
a ball $B_{1}$ of radius at least $\beta r_{0}$ inside $B_{0}\setminus A_{1}$
and so on. If at some turn Bob has no legal move or if the radii of Bob's balls do not shrink to zero we say that Alice \emph{wins by default}. Otherwise, we see
$$B_{0}\supseteq B_{0}\setminus A_{1}\supseteq B_{1}\supseteq B_{1}\setminus A_{2}\supseteq B_{2}\dots$$
and the outcome of the game is $\{x\}=\cap_{i\in\mathbb{N}}B_{i}$. A
set $S\subseteq X$ is called \emph{absolute winning} if Alice
has a strategy to make sure that either she wins by default or $x\in S$. McMullen \cite{M} showed
that in $\mathbb{R}^{d}$ absolute winning sets are winning. 
Let us consider the absolute game in $\mathcal{K}\left(\mathbb{R}^{d}\right)$.
Assume $\beta<1/3$, $r_{0}>0$ and $K\in\mathcal{K}\left(\mathbb{R}^{d}\right)$.
Then Bob can play according to the strategy that at every step he chooses a translate of
$K$ to be the center of his ball. It is easy to see that since $\beta<1/3$
no matter what Alice chooses Bob can always make a legal choice according
to his strategy (worst thing Alice can do against Bob is to choose
the same center as Bob and maximal possible radius). Hence the outcome
of the game is a translate of $K$. This  means that every absolute winning
set in $\mathcal{K}\left(\mathbb{R}^{d}\right)$ contains a translate of every
compact set. Therefore absolute winning is not the right notion
to consider when talking about generic  properties of $\mathcal{K}\left(\mathbb{R}^{d}\right)$.
\end{rem}

One of the most common applications of Schmidt games is to show existence of objects with certain properties. For example, if there are countably many    properties which are all $\alpha$-winning with a common $\alpha$, then it is also winning to satisfy all these properties at once, just like for the Banach--Mazur game.   In this paper we investigate the winning properties in terms of fractal dimensions, but many other properties could be studied. For the Banach--Mazur game there is an extensive literature studying winning properties for sets in the context of the Hausdorff metric.  However, to the best of our knowledge Schmidt's game has not yet been considered in that context. Therefore we hope to stimulate future activity in this area.

We mention one simple application of our work, mostly as motivation.  Given a notion of dimension $\dim$, one may wonder whether given any $K \subset \mathbb{R}^d$ and $s \in [0, \dim K]$ one can find a set $E \subset K$ such that $\dim E = s$.  For example, this is known to hold for upper box dimension, but \emph{not} to hold for lower box dimension, see \cite{fengcont}.  In fact, it is proved in \cite{fengcont} that there exists a compact set $K \subset [0,1]$ with the surprising property that $\underline{\dim}_B K >0$ but  that $\underline{\dim}_B E \in \{0,\underline{\dim}_B K\}$ for all $E \subset K$. It follows from our Theorem \ref{main} (i) and   Theorem \ref{extension} (ii) that, for any uniformly perfect compact set $K \subset [0,1]$, there exists a subset $E \subset K$ with $0< \dim_L E \leq  \underline{\dim}_B E<\overline{\dim}_B K$.  This proves that any set $K$ satisfying the surprising  property  demonstrated  in \cite{fengcont}  necessarily satisfies either $\dim_L K = 0$ or $ \underline{\dim}_B K<\overline{\dim}_B K$.


\subsection{Dimensions of continuous images}

Let $K$ be a compact metric space. Fix $d\in\N$, and let $\FF = \FF(K,d)$ be the space of continuous functions from $K$ to $\R^d$, endowed with the metric induced by the supremum norm which is denoted by $|\cdot|$. Write $\TD (K)$ to denote the \emph{topological dimension} of $K$.  Recall that $\TD (K)$ is the minimal integer $n \geq 0$ such that for every $r>0$ there exists an $r$-cover of $K$ by open sets such that every point $x \in K$ lies in at most $n+1$ of the covering sets. 

\begin{thm} \label{functionthm}
\hspace{1mm}
\begin{enumerate}
\item[\textup{(i)}] If $\TD(K) \geq d$, then $\{\ff\in \FF \sep \ff(K) \textup{ has nonempty interior}\}$ is winning.

\item[\textup{(ii)}] If $\TD(K) < d$, then the set $\{\ff\in \FF \sep \AD(\ff(K)) < d\}$ is winning.

\item[\textup{(iii)}] If $K$ is uncountable, then the set $\{\ff\in\FF \sep \HD(\ff(K)) > 0\}$ is winning.

\item[\textup{(iv)}] If $\TD(K) \leq d$, then the set $\{\ff\in \FF \sep \mathcal{H}^{\TD(K)}(\ff(K)) > 0\}$ is winning. In particular, the set $\{\ff\in\FF \sep \HD(\ff(K)) \geq \TD(K)\}$ is winning.
\end{enumerate}
\end{thm}

\begin{rem}
In the proof of Theorem  \ref{functionthm} we show that the conclusion of part (iii) holds for every perfect set $K$. Every uncountable compact set $K$ contains a perfect set $F$. In the more general case, Alice plays her strategy for the restriction of the functions to $F$. It follows from the general version of Tietze's extension theorem that Alice  can extend her choice to the whole set $K$ to be able to play a legal move.
\end{rem}

\begin{rem}
The conclusion of part (iv) of Theorem  \ref{functionthm} can be deduced from part (i). Alice can play her strategy in the first $\TD(K)$ coordinates completely ignoring the other coordinates, just leaving them as Bob's choice. Then they end up with a function whose image projected onto the first $\TD(K)$ coordinates is of nonempty interior. Hence the conclusion of part (iv) follows.
\end{rem}

Again, this result can be compared with the analogous results in the Banach--Mazur setting, which are due to Balka--Farkas--Fraser--Hyde \cite{balkaetal}.

\begin{thm}[Balka--Farkas--Fraser--Hyde, 2013]
The set
\[
\left\{ \ff\in \FF \sep \HD (\ff(K)) = \underline{\dim_{B}} (\ff(K))  =  \min\{ d,\, \TD(K)\} ,\, \overline{\dim_{B}} (\ff(K))=\dim_{A} (\ff (K))= d \right\}
\]
is winning for the Banach--Mazur game.
\end{thm}

In case $\TD(K) \geq d$, these properties are similar for both Schmidt winning and Banach--Mazur winning; however, in the other case they are rather different.

\begin{rem}
We can again consider the absolute game in $\mathcal{F}$. Assume $\beta<1/3$,
$r_{0}>0$ and $\ff\in\mathcal{F}$. Then Bob can play that
at every step he chooses the center of his ball on the line $\left\{\ff+\lambda\sep\lambda\in\mathbb{R}^{d}\right\}$.
It is easy to see that Bob can always make a legal choice like that
since $\beta<1/3$ (again the worst thing Alice can do against Bob is to choose
the same center as Bob and maximal possible radius). The outcome of
the game is $\ff+\lambda$ for some $\lambda\in\mathbb{R}^{d}$. This
means that an absolute winning set contains a constant translate of
every function in $\mathcal{F}$. Hence, just like in the case of
$\mathcal{K}\left(\mathbb{R}^{d}\right)$, the absolute game does not seem to
be the right notion to consider in $\mathcal{F}$.
\end{rem}

\subsection{Digit frequencies}

Finally, we consider analogous questions concerning the frequencies of digits in expansions of real numbers.  Although this is somewhat incongruous with our other results, frequencies are inherently related to densities and therefore dimension.  There are also direct connections between digit frequencies and dimension in many well-studied  settings, such as random fractals or Moran constructions, see below for a simple example.    Moreover, questions regarding the generic behaviour of digit frequencies are among the simplest and most widely studied and therefore it is useful to see how our approach fits in here. However, the most important reason that the results in this section fit with the rest of the paper is that the same phenomenon occurs: the Schmidt winning properties are starkly different from the properties which are winning for the Banach--Mazur game and, moreover, the differences are similar in spirit.

For simplicity, we consider binary expansions of numbers $x= x_0.x_1 x_2 \dots  \in \mathbb{R}$ where $x_0$ is an integer and $x_i \in \{0,1\}$ are the digits in the binary expansion of the fractional part of $x$.  For definiteness take the lexicographically maximal expansion in the situations where $x$ does not have a unique expansion and assume that all expansions are infinite, for example $1=1.000\dots$. For $x \in \mathbb{R}$ and $j \in \{0,1\}$ write
\[
d^+(x,j) = \limsup_{k \to \infty}  \frac{\#\left\{ 1 \leq i \leq k \sep x_i = j\right\}}k
\]
and $d^-(x,j)$ for the same expression with $\limsup$ replaced by $\liminf$, where if $S$ is a finite set then $\# S$ stands for the number of elements in $S$.

\begin{thm} \label{numbersthm}
The set
\[
\left\{ x \in \R \sep  0<d^-(x,j) \leq  d^+(x,j) < 1 \text{ for } j \in \{0,1\} \right\}
\]
is winning, but for all $\varepsilon>0$ and $j \in \{0,1\}$ the sets  $\left\{ x \in \R \sep  \eps <d^+(x,j) \right\} $ and
$\left\{ x \in \R \sep  d^-(x,j) < 1-\eps \right\}$ are not winning.
\end{thm}

This theorem is 
sharp in the sense that we cannot replace $0$ by any $\varepsilon>0$.

\begin{rem}
Similarly to Remark \ref{abint}, we can ask if for sets $A,B\subseteq (0,1)$, the set
\[
\left\{ x \in \R \sep d^-(x,j)\in A,\, d^+(x,j)\in B,\, d^-(x,j)<d^+(x,j)\right\}
\]
is winning if and only if $A$ and $B$ are both `big' subsets of $(0,1)$, where `big' may stand for either full Lebesgue measure, winning or co-meager.
\end{rem}

There are obvious parallels between Theorem \ref{numbersthm} and the results in the previous sections, namely Theorem \ref{main} and Corollary \ref{corcor}.  Moreover, there is a stark comparison between this result and other results concerning generic behaviour of digit frequencies.

\begin{thm}[Borel 1909, Volkmann 1959]
 The set
\[
\left\{ x \in \R \sep  d^-(x,j) = d^+(x,j) =1/2  \text{ for } j \in \{0,1\} \right\}
\]
has full Lebesgue measure and  the set
\[
\left\{ x \in \R \sep  0=d^-(x,j) \text{ and }  d^+(x,j) = 1 \text{ for } j \in \{0,1\} \right\}
\]
is winning for the Banach--Mazur game.
\end{thm}

The Lebesgue measure part of this theorem was proven by Borel in 1909 as an application of what became known as the Borel--Cantelli lemma \cite{borel} and the Banach--Mazur part was proven by Volkmann in 1959 \cite{volkmann}.

As mentioned above there is a simple (well-known) direct connection between dimension and digit frequencies. We associate a set $F \subset [0,1]$ with the binary expansion of some $x \in [0,1]$ by an iterative procedure where we begin with one level 1 interval $[0,1]$ and then at level $i \geq 1$ we replace all level $i$ intervals with two abutting  intervals of half the size, or with just one interval of half the size, according to whether the $i$th digit $x_i$ in the expansion of $x$ is 0 or 1.  Then the Hausdorff and lower box dimensions of $F$ are given by $d^-(x,0)$ and the packing and upper box dimensions of $F$ are given by $d^+(x,0)$. 

\section{Proofs}

\subsection{Proof of Theorem \ref{main}(i)}\label{sec:proofOfTheorem2part1}

Since $\dim_{L}X>0$, for all $N \geq 2$ there exists an $\alpha \in (0,1/5)$ (depending on $X$ and $N$) such that for every $y\in X$ and $0<r\leq 1$ there exist $x_1,\dots, x_{N}$ such that $\bigcup_{i=1}^{N}B(x_i,2\alpha r)\subseteq B(y,r)$
is a disjoint union. In what follows fixing $N=2$ is sufficient.

Let  $ \beta \in (0,1)$ be arbitrary. The game begins by Bob choosing a starting set $K_0^{B}\in\mathcal{K}(X)$
and a radius $r_0>0$. Alice's move is then as follows. Take $x_1,\dots,x_{M}\in K_0^{B}$, where $M \geq N$, such that $\bigcup_{i=1}^{M}B\left(x_i,2\alpha r_0\right)$
is a disjoint union, and $\bigcup_{i=1}^{M}B\left(x_i,4\alpha r_0\right)$ covers $K_0^{B}$. Alice then chooses the set $K_0^{A} =\left\{x_1,\dots, x_{M}\right\}$. This choice is legal since $d_{H}\left(K_0^{A},K_0^{B}\right)<(1-\alpha)r_0$. Then Bob
will have to choose a set $K_1^{B}$  which is contained in $\bigcup_{i=1}^{M}B\left(x_{i},(1-\beta)\alpha r_{0}\right)$
and contains at least one point in every $B\left(x_{i},(1-\beta)\alpha r_{0}\right)$
because the balls are $\alpha r_{0}$-separated. Alice then repeats
her strategy in each of the balls $B\left(x_{i},(1-\beta)\alpha r_{0}\right)$
simultaneously. At every step $n$ in each of the balls of the previous
step of the constuction Alice finds at least $N$ more new balls of
radius $(1-\beta)\alpha(\alpha\beta)^{n-1}r_{0}$. Let $K$ be the
outcome of the game, $x\in K$ and $0<R\leq r_{0}$. Let $m\in\mathbb{N}$
such that $(\alpha\beta)^{m}(1-\beta)\alpha r_{0}<R\leq(\alpha\beta)^{m-1}(1-\beta)\alpha r_{0}$.
Then $B(x,2R)$ contains at least one of the balls of step $m+1$.
Hence
\[
N_{(1-\beta)\alpha(\alpha\beta)^{n}r_{0}}(K\cap B(x,2R))\geq N^{n-m}\geq(\alpha\beta)^{-(n-m)\varepsilon}\geq\alpha\beta 2^{-1}\left(\frac{2R}{(1-\beta)\alpha(\alpha\beta)^{n}r_{0}}\right)^{\varepsilon}
\]
if $0<\varepsilon<1$ is small enough that $N\geq(\alpha\beta)^{-\varepsilon}$. In particular,
 $\dim_{L}K\geq\varepsilon>0$, which proves the result.

\subsection{Proof of Theorem \ref{main}(ii)}

Let $\varepsilon>0$ and fix $\alpha \in (0,1/2]$.  We need to show that there exists $\beta \in (0,1)$ such that
\[
\left\{ K\in\mathcal{K}(X)\sep\dim_{L}K<\dim_{L}X-\varepsilon\right\}
\]
is \emph{not} $(\alpha,\beta)$-winning. If $\dim_{L}X=0$ then the statement is trivial, so assume that $\dim_{L}X>0$. It follows that for
\[
\dim_{L}X > t > \max\left\{ \dim_{L}X-\varepsilon, 0\right\}
\]
there exist $C>0$ and $r_{0}$ such that for all $\beta \in (0,1)$, for all $0<r<r_{0}$ and all $x\in X$ there are $y_1,\dots,y_{N}$ such that $\bigcup_{i=1}^{N}B\left(y_i,2\beta r\right)\subseteq B(x,r)$ is a disjoint union, where $N\geq C\beta^{-t}$. Let $t>s>\dim_{L}X-\varepsilon$. Let $\beta>0$ be small enough that
$4\beta<1-\beta$ and $C\alpha^{s}\geq\beta^{t-s}$, i.e. $C\beta^{-t}\geq(\alpha\beta)^{-s}$.
Bob starts by choosing $K_{0}^{B}$ to be a single point and $r_{0}$ was chosen above. Then Bob adopts Alice's strategy from part (i). Similarly
it can be shown that
\[
N_{(1-\alpha)(\alpha\beta)^{n}r_{0}}(K\cap B(x,R))\geq(C\beta^{-t})^{n-m}\geq(\alpha\beta)^{-(n-m)s}\geq\alpha^{s}\beta^{s}2^{-s}\left(\frac{2R}{(1-\alpha)(\alpha\beta)^{n}r_{0}}\right)^{s}
\]
and therefore the outcome of the game satisfies $\dim_{L}K\geq s>\dim_{L}X-\varepsilon$, proving the result.

\subsection{Proof of Theorem \ref{main}(iii)}

Let $\alpha \in (0,1/2)$ and $\beta \in (0,1)$ and let $d = \dim_A X<\infty $. Since by assumption $X$ is Assouad sharp, it follows that there exists $C>0$ such that every bounded set $B$ of diameter $2R$ can be covered by fewer than
$C(R/r)^d$
balls of radius $r \in (0,R)$ which are centered in $B$, moreover we can assume that the balls with the same centers and of radius $r/2$ are disjoint. For now, assume that Bob starts the game by choosing $K_0^{B}$ to be a subset of the ball $B(y,(1-\beta)r)$ and a starting radius $r_0=\beta r$. Playing as Alice, choose $x_1,\dots,x_{N} \in K_0^{B}$ such that $K_0^{B}\subseteq\bigcup_{i=1}^{N}B(x_i,\beta r)$ where $N\leq C((1-\beta)/\beta)^{d}$ and $\bigcup_{i=1}^{N}B(x_{i},\alpha\beta r)\subseteq\bigcup_{i=1}^{N}B\left(x_{i},2^{-1}\beta r\right)$
is a disjoint union (since that $\alpha<1/2$).  Moreover,  if $\alpha<1/4$ then the distinct balls are $\alpha\beta r$-separated. Alice then chooses $K_0^A = \left\{x_1,\dots,x_{N}\right\}$. For his next move, Bob will have to choose a set $K_1^{B}$
that is contained in $\bigcup_{i=1}^{N}B\left(x_i,(1-\beta)\alpha\beta r\right)$
and the radius of the step is $\beta\alpha\beta r$. Alice in her next move repeats the previous strategy in all of the smaller balls $B(x_i,(1-\beta)\alpha\beta r)$ simultaneously and proceeds similarly at every step.

Let us denote the outcome of the game by $K$. Let $x\in K$ and $0<R\leq r_{0}$ be fixed. Let $m\in\mathbb{N}$ be
such that $r(\alpha\beta)^{m+1}<R\leq r(\alpha\beta)^{m}.$
Let's say Alice's $m$th move was $\{y_{i}\}_{i=1}^{M_m}$. Then it follows
from the strategy that $K$ is contained in $\bigcup_{i=1}^{M_m}B\left(y_{i},(1-\beta)(\alpha\beta)^{m}r\right)$
and these balls are $(\alpha\beta)^{m}r$-separated. Then the number
$M$ of these balls that intersects $B(x,R)$ is at most the number
of these balls that are contained in $B\left(x,2\alpha^{-1}\beta^{-1}R\right)$ which is at most $C\left(\frac{2\alpha^{-1}\beta^{-1}R}{(\alpha\beta)^{m}r}\right)^{d}\leq C\left(\frac{2}{\alpha \beta }\right)^{d}$
because these are disjoint balls centered in $B(x,2R)$. Let $n>m$ and observe that by the $n$th step of the game we have that each of these $M$
balls of radius $(\alpha\beta)^{m}r$ contain at most $\left(C\left(\frac{1-\beta}{\beta}\right)^{d}\right)^{n-m}$
balls of radius $(1-\beta)(\alpha\beta)^{n}r$ such that $K\cap B(x,R)$
is contained in the union of these at most $M\left(C\left(\frac{1-\beta}{\beta}\right)^{d}\right)^{n-m}\leq C\left(\frac{2}{\alpha \beta }\right)^{d}\left(C\beta^{-d}\right)^{n-m}$
balls of radius $(1-\beta)(\alpha\beta)^{n}r\leq(\alpha\beta)^{n}r$.
Hence, noting that $(\alpha\beta)^{m-n}\leq\alpha^{-1}\beta^{-2}\frac{R}{r(\alpha\beta)^{n}}$,
it follows that

\begin{align*}
N_{(\alpha\beta)^{n}r}\left(K\cap B(x,R)\right)\leq C\left(\frac{2}{\alpha \beta }\right)^{d}\left(C\beta^{-d}\right)^{n-m}
&\leq C\left(\frac{2}{\alpha \beta }\right)^{d}\left(\alpha\beta\right)^{-(d-\varepsilon)(n-m)} \\
&\leq C2^{d}\alpha^{-2d}\beta^{-3d}\left(\frac{R}{r(\alpha\beta)^{n}}\right)^{(d-\varepsilon)}
\end{align*}


provided  $0<\alpha<1/4$ is small enough that $C<\alpha^{-d/4}<\alpha^{-(d-\varepsilon)/2}$
and $0<\varepsilon<d/2$ is small enough that $\alpha^{d/2}<\left(\beta\alpha^{1/2}\right)^{\varepsilon}$.
Therefore we can choose small enough $\alpha$ such that for every $\beta$
we can choose a small enough $\varepsilon$ that $\dim_{A}K\leq d-\varepsilon$.

Finally, suppose Bob does not begin by choosing a set inside $B\left(y,(1-\beta) r\right)$, but rather an arbitrary compact set $K_0^{B}$, as he is of course permitted to do.  In this case, on her first move, Alice plays her strategy with the only exception that on the first level we have no control on the number of balls $N$. However, this is not an issue because it is just some
fixed
number and on the latter levels we have exactly the same control over the balls. It follows that for some $\alpha>0$ the set
\[
\left\{ K\in\mathcal{K}(X)\sep\dim_{A}K<\dim_{A}X\right\}
\]
is $(\alpha,\beta)$-winning for all $\beta \in (0,1)$ and therefore winning.

\subsection{Proof of Theorem \ref{main}(iv)}

Let $\varepsilon>0$ and fix $\alpha \in (0,1/2)$. We need to show that there exists $\beta \in (0,1)$ such that
\[
\left\{ K\in\mathcal{K}(X)\sep\dim_{A} K >\varepsilon\right\}
\]
is \emph{not} $(\alpha,\beta)$-winning. To achieve this, Bob will adopt a similar strategy
to Alice's strategy in the proof of Part (iii). Since $X$ is doubling, it follows that there exist $C,d>0$ such that every ball
of radius $R$ can be covered by fewer than
$C(R/r)^d$
balls of radius $r \in (0,R)$. Note that $d$ may be strictly greater then $\AD X$ in case $X$ is not Assouad sharp, but this is not an issue. Bob begins by choosing $K_0^B=\{y\}$, where $y \in X$ is arbitrary, and $r_0=1$. Alice is then forced to choose a set $K_0^A \subseteq B(y, 1-\alpha )$. Bob covers $K_{0}^{A}$ by balls of radius $\alpha r$ with centers
in $y_{1},\dots,y_{N} \in K$ where $N\leq C((1-\alpha)/\alpha)^{d}$ is such
that $\bigcup_{i=1}^{N}B(y_{i},\beta\alpha r)\subseteq\bigcup_{i=1}^{N}B\left(y_{i},2^{-1}\alpha r\right)$
is a disjoint union, provided since $\beta<1/2$.   Moreover if $\beta<1/4$
then the balls are $\beta\alpha$-separated. Alice's next choice must
be contained in $\bigcup_{i=1}^{N}B(y_{i},(1-\alpha)\beta\alpha r)$.
Bob repeats his strategy in each of these balls. Let $K$ be the outcome
of the game. Via a similar argument to that  in the proof of Part
(iii) one can deduce that for $(\alpha\beta)^{m+1}<R\leq(\alpha\beta)^{m}$
and $n>m$

\begin{align*}
N_{(\alpha\beta)^{n}}(K\cap B(x,R))\leq C\left(\frac{2}{\alpha \beta }\right)^{d}\left(C\alpha^{-d}\right)^{n-m}
&\leq C\left(\frac{2}{\alpha \beta }\right)^{d}\left(\alpha\beta\right)^{-\varepsilon(n-m)} \\
&\leq C2^{d}\alpha^{-2d}\beta^{-2d}\left(\frac{R}{(\alpha\beta)^{n}}\right)^{\varepsilon}
\end{align*}


if $0<\beta<1/4$ is small enough such  that $C\alpha^{-d}<\beta^{-\varepsilon}<\left(\alpha\beta\right)^{-\varepsilon}$. Therefore $\dim_{A}K\leq\varepsilon$, as required.

\subsection{Proof of Theorem \ref{extension} (i)}

Since $X$ is uniformly perfect, i.e. $\dim_{L}X>0$, we can find
a small enough $0<\alpha<1/4$ such that inside every ball $B$ of
radius $r/8$ we can find two disjoint balls of radius $\alpha r$
such that they are $\alpha r$-separated from each other and from
the boundary of $B$. Let $\beta\in(0,1)$ be arbitrary and choose
$\varepsilon>0$  small enough that
\begin{equation}
(\alpha\beta)^{-2\varepsilon}\leq2^{1/2}\,.\label{eq:epsikics}
\end{equation}
Write $s=\dim_{A}X<\infty$,  and therefore there exists $C_{\varepsilon}>0$
such that
\begin{equation}
N_{r}(B(x,R))\leq C_{\varepsilon}(R/r)^{s+\varepsilon}\label{eq:Nres}
\end{equation}
for every $x\in X$ and $0<r<R$. Fix $k\in\mathbb{N}$ such that
\begin{equation}
C_{\varepsilon}\leq2^{k/2}\label{eq:cepsi}
\end{equation}
and let $\rho_{n}=(\alpha\beta)^{kn}r$ where $r_{0}=\beta r$.

First assume that Bob starts by choosing $K_{0}^{B}\subseteq B(y,(1-\beta)r)\subseteq B(y,r)$
for some $y\in X$. Then Alice can find $x_{1},\dots,x_{N}\in K_{0}^{B}$
such that $K_{0}^{B}\subseteq\bigcup_{i=1}^{N}B(x_{i},\beta r/2)$ and
$\bigcup_{i=1}^{N}B(x_{i},\beta r/4)$ is a disjoint union. Then the
balls $B(x_{i},\beta r/8)$ are disjoint and $\beta r/8$-separated.
By the choice of $\alpha$, inside every ball $B(x_{i},\beta r/8)$
there exists two balls $B\left(y_{i}^{1},\alpha\beta r\right)$ and $B\left(y_{i}^{2},\alpha\beta r\right)$ which are disjoint and $\alpha\beta r$-separated from each other
and from the boundary of $B(x_{i},\beta r/8)$. Without  loss of
generality we can assume that $\widetilde{N_{\rho_{1}}}\left(B\left(y_{i}^{1},\alpha\beta r\right)\right)\leq\widetilde{N_{\rho_{1}}}\left(B\left(y_{i}^{2},\alpha\beta r\right)\right)$
where $\widetilde{N_{\rho}}(B)$ denotes the maximum number of disjoint
balls of radius $\rho$ inside $B$ for a ball $B$. Alice's move
then is $K_{1}^{A}=\bigcup_{i=1}^{N}\left\{y_{i}^{1}\right\}$.  Bob's next
move $K_{1}^{B}$ must be contained in $\bigcup_{i=1}^{N}B\left(y_{i}^{1},(1-\beta)\alpha\beta r\right)\subseteq\bigcup_{i=1}^{N}B\left(y_{i}^{1},\alpha\beta r\right)$
and the radius of the step is $\beta\alpha(\beta r)$, noting that $K_{1}^{B}$
also contains at least $1$ point from each of these balls. Assume
that on the $m$th round of the game Bob plays $K_{m}^{B}$ which, by
Alice's strategy, is contained in $\bigcup_{i=1}^{N_{m}}B\left(y_{i,m}^{1},(1-\beta )(\alpha\beta)^{m}r\right)\subseteq \bigcup_{i=1}^{N_{m}}B\left(y_{i,m}^{1},(\alpha\beta)^{m}r\right)$
where these balls are $(\alpha\beta)^{m}r$-separated. Let $m=kn+l$
where $n\geq0$, $l=0,\dots,k-1$. Then in each of these balls similarly
to the first step of the game Alice can find $x_{1}^{m,i},\dots,x_{N_{m}^{i}}^{m,i}\in K_{0}^{B}\cap B\left(y_{i,m}^{1},(\alpha\beta)^{m}r\right)$
such that $K_{0}^{B}\cap B\left(y_{i,m}^{1},(\alpha\beta)^{m}r\right)\subseteq\bigcup_{i=1}^{N}B\left(x_{j}^{m,i},\beta(\alpha\beta)^{m}r/2\right)$
and $\bigcup_{i=1}^{N_{m}^{i}}B\left(x_{j}^{m,i},\beta(\alpha\beta)^{m}r/4\right)$
is a disjoint union. Then the balls $B\left(x_{j}^{m,i},\beta(\alpha\beta)^{m}r/8\right)$
are disjoint and $\beta(\alpha\beta)^{m}r/8$-separated. By the choice
of $\alpha$ inside every ball $B\left(x_{j}^{m,i},\beta(\alpha\beta)^{m}r/8\right)$
there exist  two balls $B\left(y_{j}^{1,m,i},(\alpha\beta)^{m+1}r\right)$ and
$B\left(y_{j}^{2,m,i},(\alpha\beta)^{m+1}r\right)$ which are disjoint and
$(\alpha\beta)^{m+1}r$-separated from each other and from the boundary
of $B(x_{i},\beta r/8)$. Without  loss of generality we can assume
that
\begin{equation}
\widetilde{N_{\rho_{n+1}}}\left(B\left(y_{j}^{1,m,i},(\alpha\beta)^{m+1}r\right)\right)\leq\widetilde{N_{\rho_{n+1}}}\left(B\left(y_{j}^{2,m,i},(\alpha\beta)^{m+1}r\right)\right).\label{eq:dijointingnr}
\end{equation}
Alice's move is then $K_{m+1}^{A}=\bigcup_{i=1}^{N_{m}}\bigcup_{j=1}^{N_{m}^{i}}\left\{y_{j}^{1,m,i}\right\}$=$\bigcup_{i=1}^{N_{m+1}}\left\{y_{i,m+1}^{1}\right\}$
and Bob's next move must be contained in $\bigcup_{i=1}^{N_{m+1}}B\left(y_{i,m+1}^{1},(1-\beta)(\alpha\beta)^{m+1}r\right)$.

Now let us consider what happens inside one ball of the $nk$th level
of the construction $B=B\left(y_{i,nk}^{1},(\alpha\beta)^{nk}r\right)$ at the
$(n+1)k$th level of the construction. The ball $B\left(y_{i,nk}^{1},(\alpha\beta)^{nk}r\right)$
contains some $M$ balls from the collection $B\left(y_{i,(n+1)k}^{1},(\alpha\beta)^{(n+1)k}r\right)$.
On the $((n+1)k-1)$th level we kept the balls $B\left(y_{j}^{1,(n+1)k-1},(\alpha\beta)^{(n+1)k-1}r\right)$
and eliminated the balls $B\left(y_{j}^{2,(n+1)k-1},(\alpha\beta)^{(n+1)k-1}r\right)$
and by \eqref{eq:dijointingnr} it follows that inside the balls at the $((n+1)k-1)$th level that are contained in $B$ there are at least $2M$
balls of radius $(\alpha\beta)^{(n+1)k}r$. Similarly inside the balls
at the $((n+1)k-2)$the level of the construction that are contained
in $B$ there are at least $4M$ balls of radius $(\alpha\beta)^{(n+1)k}r$.
Continuing this through $k$ steps we get that inside $B$ there are
$2^{k}M$ balls of radius $(\alpha\beta)^{(n+1)k}r$. Hence
\begin{equation}
M\leq2^{-k}N_{(\alpha\beta)^{(n+1)k}r}(B)\leq2^{-k}C_{\varepsilon}(R/r)^{s+\varepsilon}=2^{-k}C_{\varepsilon}(\alpha\beta)^{-2\varepsilon k}(\alpha\beta)^{-(s-\varepsilon)k}\leq(\alpha\beta)^{-(s-\varepsilon)k}\label{eq:Mupper}
\end{equation}
by \eqref{eq:epsikics}, \eqref{eq:Nres} and \eqref{eq:cepsi}. So
inside each ball at the $nk$th level of the construction there are
at most $(\alpha\beta)^{-(s-\varepsilon)k}$ balls from the $(n+1)k$
level of the construction while the ratio between the radii
are $(\alpha\beta)^{k}$. Hence for the outcome of the game $K$
\[
\dim_{A}K\leq s-\varepsilon\,.
\]
Now let us consider the situation when Bob does not start by choosing
$K_{0}^{B}\subseteq B(y,(1-\beta)r)$. Then Alice plays her strategy
exactly the same way. The only difference is that now on the first
level, and so on the first $k$ levels, we have no control over the number
of balls of the construction. However, from then on inside each ball
we have exactly the same control over the number of the balls. Hence
this has no effect on $\dim_{A}K$ and the proof is finished. Note
that
\begin{equation}
N_{\rho_{n}}(K)\leq N_{(\alpha\beta)^{nk}r}\left(B\left(K_{nk}^{A},\rho_{n}\right)\right)\leq(\alpha\beta)^{-(s-\varepsilon)nk}N_{\beta r/2}\left(K_{0}^{B}\right)=\rho_{n}^{-(s-\varepsilon)}N_{\beta r/2}\left(K_{0}^{B}\right).\label{eq:laterprofeq}
\end{equation}

\subsection{Proof of Theorem \ref{extension} (ii)}

We prove (ii) similarly to (i). We apply a trick to replace $\dim_{A}X$
with $\overline{\dim_{B}}X$, but the price is that we cannot keep
$k$ fixed and therefore only get estimates for the lower box dimension
of the outcome. The Assouad spectrum, introduced by Fraser and Yu
$\cite{Spectra}$, of a compact set $X$ is the function $\theta\mapsto\dim_{A}^{\theta}X$,
where $\theta$ varies in $(0,1)$ and
\begin{multline*}
\dim_{A}^{\theta}X := \inf\bigg\{\alpha \sep \text{there exists }C>0\text{ such that all }0<R<1\text{ and }x\in X\text{ satisfy } \\
N_{R^{1/\theta}}\left(B(x,R)\right)\ \leq\ C\left(\frac{R}{R^{1/\theta}}\right)^{\alpha}\bigg\}\,.
\end{multline*}
It was shown in $\cite{Spectra}$ that $\dim_{A}^{\theta}X$ is continuous
in $\theta$ and satisfies
\[
\dim_{A}^{\theta}X\leq\frac{\overline{\dim_{B}}X}{1-\theta}
\]
and therefore $\dim_{A}^{\theta}X\to\overline{\dim_{B}}X$ as $\theta\to0$.
Therefore, in the above proof we can replace \eqref{eq:Nres} with
the statement that there exist a small $\theta\in(0,1)$ and a constant
$C_{\varepsilon,\theta}>0$ such that
\begin{equation}
N_{R^{1/\theta}}\left(B(z,R)\right)\leq C_{\varepsilon,\theta}\left(R/R^{1/\theta}\right)^{s+\varepsilon}\label{assfact2}
\end{equation}
for all $z\in X$ and $R\in(0,1)$ where $s$ is now the upper box
dimension of $X$, that is, $s=\overline{\dim_{B}}X$ (note that we
choose $\theta$ to be less than $1/2$ which ensures that $1/\theta-1\geq1$).
The only place this estimate is needed in the proof of (i) is \eqref{eq:Mupper}
but to make $\eqref{assfact2}$ work here we would ideally choose
\[
\rho_{n+1}=\rho_{n}^{1/\theta}
\]
which means that $k$ needs to change at every step of the strategy.
That is, at step $n$ we group the turns in Schmidt's game into groups
of length $k_{n}$ defined inductively as follows. Let $k\in\mathbb{N}$
be such that
\[
C_{\varepsilon,\theta}\leq2^{k/2}\,.
\]
We need to make sure that $k_{n}\geq k$ for every $n$. Choose $k_{1}=k$ and,
given $k_{1},\dots,k_{n}\geq k$, choose $k_{n+1}$ to be the largest
integer satisfying
\[
\rho_{n+1}:=(\alpha\beta)^{\sum_{i=1}^{n+1}k_{i}}r\geq\left((\alpha\beta)^{\sum_{i=1}^{n}k_{i}}r\right)^{1/\theta}=\rho_{n}^{1/\theta}\,.
\]
It is easy to see that if $r\leq1$ (which can be assumed because
otherwise Alice just plays randomly until $r_{m}/\beta\leq1$ and
starts playing her strategy then), then
\[
k_{n+1}\geq(1/\theta-1)\sum_{i=1}^{n}k_{i}\geq nk\geq k\,.
\]
Clearly this forces $k_{n}$ to grow very quickly and, moreover,
\[
\rho_{n}^{1/\theta}(\alpha\beta)^{-1}\geq\rho_{n+1}\geq\rho_{n}^{1/\theta}\,.
\]
This is enough to apply $\eqref{assfact2}$ to make an estimation
like in \eqref{eq:Mupper} for all $n$ up to a constant factor which
can be absorbed into the constant $C_{\varepsilon,\theta}$.

The proof then works as before leading to the estimate (similar to
\eqref{eq:laterprofeq})
\[
N_{(\alpha\beta)^{S(n)}r}\left(K\right)\leq(\alpha\beta)^{-S(n)(s-\varepsilon)}N_{\beta r/2}\left(K_{0}^{B}\right)
\]
for all $n\in\mathbb{N}$ where $S(n)=\sum_{i=1}^{n}k_{i}$. Therefore
we conclude that $\underline{\dim_{B}}K\leq s-\varepsilon<\overline{\dim_{B}}X$
completing the proof.


\subsection{Proof of Theorem \ref{functionthm}(i)}

The following notation is used in this section: vectors in $\bbr^d$ and functions to $\bbr^d$ are denoted by boldface letters and their coordinates by a normal font with a subscript, e.g., for $\ff:K\to\bbr^d$, the functions $f_1,\ldots,f_d:K\to\bbr$ are such that $\ff=\left(f_1,\ldots,f_d\right)$.

Fix any $\alpha\in(0,1/3)$ and $\beta\in(0,1)$, and let $\ff\in\FF$ and $r_0 > 0$ be Bob's first move and starting radius. Write $r = (1-\alpha)r_0$, and for each $\bfm\in\Z^d$ let
\[
K_{\bfm} = \bigcap_{i=1}^df_{i}^{-1}\left([m_i r, (m_i+1) r]\right).
\]
Then $K = \bigcup_{\bfm\in\Z^d} K_{\bfm}$, so by the sum theorem \cite[Theorem 1.5.3]{Engelking}, there exists $\bfm\in\Z^d$ for which
\[
\TD(K_{\bfm}) \geq d\,.
\]
Therefore, by \cite[Corollary of Theorem II.8]{nagata}, there exist closed sets $F_1,\ldots,F_d \subseteq K_{\bfm}$ and relatively open sets $U_1,\ldots,U_d \subseteq K_{\bfm}$ such that $F_i \subseteq U_i$ for each $1\leq i\leq d$, and such that if $F_i \subseteq W_i \subseteq U_i$ for some relatively open sets $W_1,\ldots,W_d \subseteq K_{\bfm}$, then $\bigcap_1^d \del W_i \neq \emptyset$.
Now Alice's strategy is as follows: on the first move,
choose $\gg\in B\left(\ff,(1-\alpha)r_0\right)$ so that $g_i = m_i r$ on $F_i$ and $g_i = (m_i+1) r$ on $K_{\bfm}\butnot U_i$ for some $\bfm\in\Z^d$ for every $1\leq i \leq d$, and play arbitrarily on later moves. The outcome function $\hh$ must lie in $B\left(\gg,\alpha r_0\right)$, hence it satisfies
\begin{align*}
h_i &\leq \left(m_i + \frac{\alpha}{1-\alpha}\right)r && \text{  on  } F_i\,,\\
h_i &\geq \left(m_i + \frac{1-2\alpha}{1-\alpha}\right)r && \text{  on  } K_{\bfm}\butnot U_i\,,
\end{align*}
for every $1\leq i \leq d$. It follows that $\hh\left(K_{\bfm}\right)$ has a nonempty interior as it contains the cube
\[
P := \prod_{i=1}^d\left(\left(m_i + \frac{\alpha}{1-\alpha}\right)r, \left(m_i + \frac{1-2\alpha}{1-\alpha}\right)r\right).
\]
Indeed, if $\pp\in P$, then for each $i$ the set $W_i = \{h_i < p_i\}$ is open and satisfies $F_i \subseteq W_i \subseteq U_i$. Thus by our hypothesis there exists $x \in \bigcap_1^d \del W_i$, and it follows that $\hh(x) = \pp$.

\subsection{Proof of Theorem \ref{functionthm}(ii)}

Fix any $\alpha\in(0,1/4)$ and $\beta\in(0,1)$, and let $r_0 > 0$ be Bob's initial radius. Let $S \subseteq \R^d$ be the set of all points with at least one integer coordinate, i.e.
$$S=\bigcup_{i=1}^d\left\{\bx\in\bbr^d \sep x_i\in \mathbb{Z}\right\}.$$
Alice's strategy is as follows: On turn $n$, choose a ball with center $\ff$ such that $\ff(K) \subseteq \frac{r_n}{2} S$.
Such a choice is possible because $\TD(K) < d$. 
Indeed, $\TD(K) < d$ implies that the set $\{\ff\in\FF \sep \ff(K) \text{ is nowhere dense}\}$ is dense in $\FF$. 
This follows from \cite[Theorem 2.4]{balkaetal} since if $\ff(K)$ is not nowhere dense, then $\BD \ff(K) = d$.  If Bob's $n$th move is a ball with center $\gg$, let $\hh\in\FF$ be such that $\hh(K)$ is nowhere dense and $|\hh-\gg| < \left(1/2-\alpha\right)r_n$. Let $\PP$ be the collection of cubes whose boundaries form $\frac{r_n}{2} S$. For each $P\in\PP$, choose $\xx\in P\butnot \hh(K)$, and let $\pi_P:P\butnot\{\xx\}\to\del P$ be the radial projection from $\xx$. Finally, let $\ff = \pi\circ\hh$, where $\pi(\yy) = \pi_P(\yy)$ for all $\yy\in P$. Alice can choose a ball with center $\ff$ as her next move since
\[
|\ff - \gg| \leq |\ff - \hh| + |\hh - \gg| \leq (1-\alpha)r_n\,.
\]
Using this strategy Alice guarantees that the outcome of the game is a function whose image is a subset of
\begin{equation}\label{eq:porousSet}
F_d = \bigcap_{n=1}^\infty B\left(\frac{r_n}{2} S,\alpha r_n\right).
\end{equation}
The set $F_d$ is \emph{porous}, i.e., there exists $\varepsilon>0$ such that for all
$x \in \bbr^d$ and $r>0$,
there exists $y\in\bbr^d$ such that
$B(y, \varepsilon r) \subseteq B(x,r)$ and $B(y, \varepsilon r) \cap F_{d} = \emptyset$.
First we show that $F_1$ is porous. Set
$$\eps=\alpha\beta(1/2-2\alpha)$$
and assume $r>0$. Fix $n\geq1$ to be the unique integer that satisfies $r_n \leq r <r_{n-1}$. By \eqref{eq:porousSet}, $F_1$ is contained in $B\left(\frac{r_n}{2} \bbz,\alpha r_n\right)$
which is a union of intervals of length $\alpha r_n$ for which the centers are at a distance of $r_n/2$ apart. Hence, the length of the gap between consecutive intervals is $(1/2-2\alpha)r_n$. Since $r\geq r_n$, any ball of radius $r$ contains such a gap. Since $r<r_{n-1}$ we get that
\[
\frac{(1/2-2\alpha)r_n}{r} \geq \frac{(1/2-2\alpha)r_n}{r_{n-1}} = \eps\,,
\]
which finishes the proof of the one-dimensional case.
So we showed that for the scale $r$ there is a hole of size $\varepsilon r$ in $F_1 \cap B(x,r)$. If we take the product of $d$ holes of $F_1$ then it forms a hole in $F_d$. Hence in every ball we can find a hole of comparable size (the constant $\varepsilon$ depends on $\alpha,\beta$ and $d$). This shows the porosity of $F_d$.

It follows from \cite[Theorem 2.4]{fraseryu} that $\AD(F_d) < d$, so the image of the outcome has Assouad dimension strictly smaller than $d$.   Specifically, it follows from \cite[Theorem 2.4]{fraseryu} that if the Assouad dimension of a set $E \subset \mathbb{R}^d$ is $d$, then the unit ball is a \emph{weak tangent} to $E$. Having the unit ball as a weak tangent clearly prevents the set $E$ from being  porous.


\subsection{Proof of Theorem \ref{functionthm}(iii)}

Fix $\alpha<1/8$ small enough and let $S\subseteq\bbr^d$ be a $5\alpha$-separated subset of the Euclidean ball of radius $1/4$ around zero in $\bbr^d$. Set $N = \# S$. We can assume that   $N\geq2$ by choosing  $\alpha$ small enough. Alice's strategy will involve constructing a sequence of sets $E_n \subseteq K$, such that $\# E_n = N^n$. The set $E_0$ is any singleton. Suppose that the set $E_n$ has been constructed, and that the game has been played up until turn $n$, with $\gg\in\FF$ being the center of the current ball, i.e. Bob's last move. Also suppose that
\begin{equation}\label{indhyp}
|\gg(y) - \gg(x)| > 2\alpha r_{n-1} \text{ for all }x,y\in E_n \text{ distinct}\,.
\end{equation}
For each $x\in E_n$, let $F_{n,x}\subseteq K$ be a set of cardinality $N$ such that for all $y\in F_{n,x}$,
\[
|\gg(y) - \gg(x)| \leq r_n/4\,.
\]
Such a set exists because $K$ is perfect and $\gg$ is continuous. By Tietze's extension theorem there exists $\gg_{0}\in\FF$ such that $\gg_{0}(y)=\gg(x)$ for every $y\in F_{n,x}$ and $x\in E_n$ and
\[
|\gg(y) - \gg_{0}(y)| \leq r_n/4
\]
for every $y\in K$. Now let $E_{n+1} = \bigcup_{x\in E_n} F_{n,x}$. By \eqref{indhyp}, the sets $\left\{F_{n,x}\right\}_{x\in E_n}$ are disjoint, so $\# E_{k+1} = N^{k+1}$. Alice will choose the center $\ff\in\FF$ of her next ball so as to guarantee that
\begin{equation}\label{eq:confusing}
\{\ff(y) \sep y\in F_{n,x}\} = \gg(x) + r_n S
\end{equation}
for all $x\in E_n$. Namely, for each $x\in E_n$, let $\sigma_x:F_{n,x}\to S$ be any bijection, let $\left\{\phi_y\right\}_{y\in E_{n+1}}$ be any family of continuous real valued functions with disjoint supports such that $|\phi_y|\leq 1$ and $\phi_y(y) = 1$. Then the function
\begin{equation}\label{eq:confusing2}
\ff := \gg_{0} + r_n\sum_{x\in E_n} \sum_{y\in F_{n,x}} \phi_y \sigma_x(y)
\end{equation}
satisfies \eqref{eq:confusing}. Now \eqref{eq:confusing2} implies that
\[
|\ff-\gg|\leq|\ff-\gg_{0}|+|\gg_{0}-\gg|\leq\frac{r_n}{4}+\frac{r_n}{4}\leq(1-\alpha)r_n
\]
so $\ff$ is a legal move. If $\hh$ is the center of Bob's next ball, then $\hh$ must satisfy $|\hh-\gg|<(1-\beta)\alpha r_n<\alpha r_n$. Since $S$ is $5\alpha$-separated it follows that
\[
|\hh(y) - \hh(x)| \geq |\ff(y) - \ff(x)| - 2\alpha r_n > 2\alpha r_n
\]
for all $x,y\in F_{n,z}$ distinct and $z\in E_{n}$. If $x\in F_{n,z_1}$ and $y\in F_{n,z_2}$ for distinct  $z_1,z_2\in E_{n}$ then
\[
|\hh(y) - \hh(x)| \geq |\ff(y) - \ff(x)| - 2\alpha r_n \geq \alpha r_{n-1}-2r_n/4-2\alpha r_n
\]
\[
 \geq(1-1/2-2\alpha )\alpha r_{n-1}\geq \alpha r_{n-1}/4>2\alpha r_n
\]
by \eqref{indhyp}, \eqref{eq:confusing} and that $\alpha <1/8$. Thus \eqref{indhyp} is satisfied for $n+1$, allowing Alice to
complete the game by induction.

Alice's strategy results in a branching construction. Namely, if $\gg_{n}$
is Bob's $n$th move and $E_{n}$ has been constructed then the collection
of balls in the branching is $\bigcup_{x\in E_{n}}B(\gg(x),r_{n})$.
The balls in the construction are clearly disjoint by  \eqref{indhyp}. Then Alice
chooses $\ff$ such that $\ff(F_{n,x})\subseteq B(g(x),r_{n}/4)$
and for his next move Bob can only change the values by $(1-\beta)\alpha r_{n}$
hence the values stay inside $B(\gg(x),r_{n}/2)$. Thus the balls of
radius $r_{n+1}$ around these modified values are contained in $B(\gg(x),r_{n})$.
So indeed $\bigcup_{x\in E_{n}}B(\gg(x),r_{n})$ is a branching construction
with limit set contained in the image of the outcome of the game.
Hence by the mass distribution principle the image of the outcome
of the game has Hausdorff dimension at least
\[
\frac{\log N}{-\log\alpha\beta}>0\,.
\]

\subsection{Proof of Theorem \ref{numbersthm} }

We first prove that the set
\[
X: = \left\{ x \in \R \sep  0<d^-(x,0) \right\}
\]
is winning. Fix $\alpha \in (0,1/8), \beta \in (0,1)$. Bob starts the game by choosing $x_0 \in \mathbb{R}$ and $r_0>0$. As usual, we can assume $r_0 < 1/2$ since if it was not, Alice could play arbitrarily until $r_n < 1/2$. Let $I_0$ be a dyadic interval at level $k_0: = \lceil -\log r_0/\log 2 \rceil$ which is completely contained inside $B(x_0,r_0)$. Let
\[
l_0 = \max\left\{k \in \mathbb{N} \sep 2^{-k_0-k}> 2\alpha r_0\right\}
\]
and let $y_0$ be the centre of the leftmost level $(k_0+l_0)$ dyadic interval inside $I_0$. In particular, $B(y_0, \alpha r_0) \subseteq I_0 \subseteq B(x_0,r_0)$ and for all numbers in $B(y_0, \alpha r_0)$ all of the binary digits from the $(k_0+1)$th position to the $(k_0+l_0)$th position are equal to 0.

This argument is then repeated. Assuming Bob has chosen $x_n$ as the center of his $n$th ball, let $I_n$ be a dyadic interval at level $k_n: = \lceil -\log r_n/\log 2 \rceil $ which is completely contained inside $B(x_n,r_n)$. Let
\[
l_n = \max\left\{k \in \mathbb{N} \sep 2^{-k_n-k}> 2\alpha r_n\right\}
\]
and let $y_n$ be the centre of the leftmost level $(k_n+l_n)$ dyadic interval inside $I_n$. In particular, $B(y_n, \alpha r_n) \subseteq I_n \subseteq B(x_n,r_n)$ and for all numbers in $B(y_n, \alpha r_n)$ all of the binary digits from the $(k_n+1)$th position to the $(k_n+l_n)$th position are equal to 0.

Let $x$ be the outcome of the game and note that $l_n \geq -\log(8\alpha)/\log 2 >0$.  It follows that
\begin{align*}
d^-(x,0) = \liminf_{k \to \infty}  \frac{\#\left\{ 1 \leq i \leq k \sep x_i = 0\right\}}k
& \geq  \liminf_{n \to \infty}  \frac{\#\left\{ 1 \leq i \leq k_n \sep x_i = 0\right\}}{k_{n+1}} \\
& \geq  \liminf_{n \to \infty}   \frac{\sum_{i=0}^{n-1}l_i }{k_{n+1}} \\
& \geq  \liminf_{n \to \infty}   \frac{ -n \log(8\alpha)/\log 2 }{1-\log ((\alpha \beta)^{n+1} r_0) /\log 2} \\
& =     \frac{\log(8 \alpha)}{\log(\alpha\beta)} >0
\end{align*}
which proves that $X$ is $(\alpha, \beta)$-winning and therefore winning. Proving that the other sets from the theorem are winning is very similar and so we omit the proofs.


Let $\eps>0$. We now prove that the set
\[
Y: = \left\{ x \in \R \sep  \eps <d^+(x,0) \right\}
\]
is not winning. However, this is also very similar to the above proof, but with Bob adopting Alice's strategy and therefore we only point out the differences.

Fix $\alpha \in (0,1/2)$ and choose $\beta \in (0,1 )$ such that
\[
1-\frac{\log 8\beta}{-\log(\alpha \beta)} \leq \eps\,.
\]
Bob starts the game by choosing $x_0=1 \in \mathbb{R}$ and $r_0 =1$. Alice then chooses 
$y_1 \in \mathbb{R}$ such that $B(y_1, \alpha ) \subseteq B(x_0, 1)$. Let $I_1'$ be a dyadic interval at level $k_1': = \lceil -\log \alpha/\log 2 \rceil$ which is completely contained inside $B(y_1,\alpha )$. Let
\[
l_1'= \max\left\{k \in \mathbb{N} \sep 2^{-k_1'-k}> 2\alpha \beta  \right\}
\]
and Bob then chooses $x_1$ to be the centre of the rightmost level $(k_1'+l_1')$ dyadic interval inside $I_1'$. In particular, $B(x_1, \alpha \beta ) \subseteq I_1' \subseteq B(y_1, \alpha )$ and for all numbers in $B(x_1, \alpha \beta )$, all of the binary digits from the $(k_1'+1)$th position to the $(k_1'+l_1')$th position are equal to 1.

This argument is then repeated. When Alice chooses $y_n$ as the center of her $n$th ball, let $I_n'$ be a dyadic interval of level $k_n': = \lceil -\log (\alpha r_n)/\log 2 \rceil$ which is completely contained inside $B(y_n, \alpha r_n)$. Let
\[
l_n' = \max\left\{k \in \mathbb{N} \sep 2^{-k_n'-k}> 2\alpha \beta r_n\right\}
\]
and let 
$x_n$ be the centre of the rightmost level $(k_n'+l_n')$ dyadic interval inside $I_n'$. In particular, $B(x_n, \alpha \beta r_n) \subseteq I_n' \subseteq B(y_n,\alpha r_n)$ and for all numbers in $B(x_n, \alpha\beta  r_n)$ all of the binary digits from the $(k_n'+1)$th position to the $(k_n'+l_n')$th position are equal to 1.

Let $x$ be the outcome of the game and note that $l_n'  \geq -\log(8\beta )/\log 2 >0$.  It follows that
\begin{align*}
d^+(x,0) = \limsup_{k \to \infty}  \frac{\#\left\{ 1 \leq i \leq k \sep x_i = 0\right\}}k
&\leq  \limsup_{n \to \infty}  \frac{\#\left\{ 1 \leq i \leq k_n' \sep x_i = 0\right\}}{k_{n-1}'} \\
&\leq  \limsup_{n \to \infty}  \frac{k_n' - \sum_{i=0}^{n-1}l_n'}{k_{n-1}'} \\
&\leq  1- \liminf_{n \to \infty}  \frac{ -n \log(8\beta )/\log 2 }{-\log (\alpha (\alpha \beta)^{n-1}) /\log 2} \\
&= 1-    \frac{\log 8\beta}{-\log(\alpha \beta)}  \leq \eps
\end{align*}
which proves that $Y$ is not $(\alpha, \beta)$-winning and therefore not winning.

Proving that the other sets from the theorem are not winning is very similar and so we omit the proofs.

We note that Alice's and Bob's strategy in the above proofs ensure that the outcome is as desired for at least one valid binary expansion. The proof so far did not pay attention to the fact that a number may have two valid binary expansions: if eventually there are only $1$s in the binary expansion then it has another expansion  where it has only $0$s eventually. However, it is easy to modify the proofs above so that no outcome has two different valid expansions.  For example,  every now and then (but rarely) the players play a round with the opposite digit than they usually play for. We exclude the exact details.

\section*{Acknowledgements}

We would like to thank the anonymous referee for the thorough readthrough and useful comments. This work began during the semester programme on \emph{Fractal Geometry, Hyperbolic Dynamics and Thermodynamical Formalism} hosted by ICERM in Spring 2016.  It continued at the  semester programme on \emph{Fractal Geometry and Dynamics} hosted by the  Institut Mittag-Leffler in Fall 2017.  AF was financially supported by an \emph{ERC Consolidator Grant} (772466) and by \emph{The MTA Momentum Project} (LP2016-5).  JMF was financially supported by a  \emph{Leverhulme Trust Research Fellowship} (RF-2016-500) and  an \emph{EPSRC Standard Grant} (EP/R015104/1). EN and DS were supported by an  EPSRC Programme Grant (EP/J018260/1).
\\

{\small
\'Abel Farkas, E-mail: \texttt{thesecondabel@gmail.com}
\vskip.2em
Jonathan M. Fraser, E-mail: \texttt{jmf32@st-andrews.ac.uk}
\vskip.2em
Erez Nesharim, E-mail: \texttt{ereznesh@gmail.com}
\vskip.2em
David Simmons, E-mail: \texttt{david9550@gmail.com}
\vskip.2em
}

\end{document}